\documentclass[a4paper, 11pt]{article}
\usepackage[english]{babel}
\usepackage{latexsym}
\usepackage{amsfonts, amsthm}
\usepackage{amssymb}
\usepackage{amscd}
\usepackage{amsmath}
\usepackage{graphics}
\usepackage{eepic}
\usepackage{epsfig}
\usepackage{xtheorem}
\usepackage{enumerate}
\newtheorem{plain}{theorem}{Theorem}[section]
\newtheorem{definition}{definition}[theorem]{Definition}
\newtheorem{plain}{proposition}[theorem]{Proposition}
\newtheorem{plain}{lemma}[theorem]{Lemma}
{Corollary}
\newtheorem{remark}{remark}[theorem]{Remark}
{Example}

\newcommand{\C}{{\mathbb C}}

\begin{document}

\title{\textbf{\large{POLES OF THE TOPOLOGICAL ZETA FUNCTION ASSOCIATED TO AN IDEAL IN DIMENSION TWO}}}
\author{Lise Van Proeyen and Willem Veys
\footnote{K.U.Leuven, Departement Wiskunde, Celestijnenlaan 200B,
B-3001 Leuven, Belgium, email: Lise.VanProeyen@wis.kuleuven.be,
Wim.Veys@wis.kuleuven.be.  The research was partially supported by
the Fund of Scientific Research - Flanders (G.0318.06). The
original publication is available at www.springerlink.com .} }
\date{}

\maketitle

\begin{abstract}
To an ideal in $\mathbb{C}[x,y]$ one can associate a topological
zeta function. This is an extension of the topological zeta
function associated to one polynomial. But in this case we use a
principalization of the ideal instead of an embedded resolution of
the curve.
\\
\indent In this paper we will study two questions about the poles
of this zeta function. First, we will give a criterion to
determine whether or not a candidate pole is a pole. It turns out
that we can know this immediately by looking at the intersection
diagram of the principalization, together with the numerical data
of the exceptional curves. Afterwards we will completely describe
the set of rational numbers that can occur as poles of a
topological zeta function associated to an ideal in dimension two.
The same results are valid for related zeta functions, as for
instance the motivic zeta function.
\end{abstract}

 \begin{center} \footnotesize{\emph{2000 Mathematics Subject Classification.} 14E15,
14H20, 32S05.}
\end{center}

\section{Introduction}

We will first define the topological zeta function for one
polynomial in $n$ variables over $\mathbb{C}$ and mention a number
of important results about the poles of these functions.
Afterwards, we will concentrate on the topological zeta function
associated to an ideal in $\mathbb{C}[x,y]$ and make some similar
statements about its poles.

\vspace{0.5cm}

Let $f \in \mathbb{C}[x_1, \ldots, x_n]$ be a non-constant
polynomial satisfying $f(0)=0.$ To define the topological zeta
function $Z_{top, f}(s),$ we take an embedded resolution $h : X
\to \mathbb{C}^n$ of $f^{-1}\{0\}.$ Let $E_i$ for $i \in S$ be the
irreducible components of $h^{-1}(f^{-1}\{0\}),$ then we denote by
$N_i$ and $\nu_i-1$ the multiplicities of $E_i$ in the divisor on
$X$ of $f \circ h$ and $h^*(dx_1 \wedge \ldots \wedge dx_n),$
respectively. (Further on we give a description of these
multiplicities with local coordinates.) With these numerical data
we can define the local \emph{topological zeta function associated
to $f:$}
$$Z_{top,f}(s):= \sum_{I \subset S} \chi(E_I^\circ \cap
h^{-1}\{0\}) \prod_{i \in I} \frac1{N_is+\nu_i},$$ where
$\chi(\cdot)$ denotes the topological Euler-Poincar\'e
characteristic and $E_I^\circ := (\cap_{i \in I}E_i) \backslash
(\cup_{j \not \in I} E_j).$

There is also a \emph{global} topological zeta function, where we
replace $E_I^\circ \cap h^{-1}\{0\}$ by $E_I^\circ.$ Denef and
Loeser proved in \cite{DL-zeta fie onafh van res} that these
definitions are independent of the choice of the resolution.

In particular, the poles of the topological zeta function of $f$
are interesting numerical invariants. For example, the monodromy
conjecture relates the poles with eigenvalues of the local
monodromy of $f$ (see e.g. \cite{DL-zeta fie onafh van res}). It
is easy to see that all poles belong to the set $\{ -\nu_i / N_i
\, | \, i \in S \}.$ These elements are called the candidate poles
associated to the given resolution. They are all negative rational
numbers. It is an important question to determine whether or not a
candidate pole is a pole.

In \cite{veys polen top zetafunctie}, the second author proved a
fast criterion to answer this question if we work with a curve $f
\in \mathbb{C}[x_1, x_2]$. He showed that we can read the poles
out of the minimal embedded resolution of the curve: a candidate
pole $s_0$ is a pole if and only if $s_0 = -\frac{\nu_i}{N_i}$ for
some exceptional curve $E_i$ intersecting at least three times
other components or $s_0 = -\frac{1}{N_i}$ for some irreducible
component $E_i$ of the strict transform of $f.$

There are also various results about the set $$\mathcal{P}_n :=
\{s_0 \, | \, \exists f \in \mathbb{C}[x_1, \ldots , x_n] :
Z_{top,f}(s) \mbox{ has a pole in } s_0\}.$$ For example, in
\cite{artikel Ann & co} it is shown that each rational number in
the interval $[-(n-1)/2,0)$ is contained in $\mathcal{P}_n.$ For
$n=2$ this means that we know $\mathcal{P}_2$ completely, as in
\cite{SV smallest poles} it is proven that $\mathcal{P}_2 \,  \cap
\ (-\infty , -1/2\, ) \  = \{-1/2-1/i \, | \, i \in
\mathbb{Z}_{>1} \}.$

\vspace{0.5cm}

The construction of blowing-up that is used to desingularize
varieties, can also be used to principalize an ideal. This means
that after these blow-ups, the ideal is locally principal and
monomial. This is a result of Hironaka \cite{Hironaka}.

\begin{theorem} \emph{(Hironaka.)}
Let $X_0$ be a smooth algebraic variety over a field of
characteristic zero, and $\mathcal{I}$ a sheaf of ideals on $X_0.$
There exists a principalization of $\mathcal{I},$ that is a
sequence
$$X_0 \stackrel{\sigma_1}{\longleftarrow} X_1 \stackrel{\sigma_2}{\longleftarrow} X_2
\cdots \stackrel{\sigma_i}{\longleftarrow} X_i \longleftarrow
\cdots \stackrel{\sigma_r}{\longleftarrow} X_r=X$$ of blow-ups
$\sigma_i: X_{i-1} \leftarrow X_i$ in smooth centers $C_{i-1}
\subset X_{i-1}$ such that
\begin{enumerate}
\item the exceptional divisor $E^i$ of the induced morphism
$\sigma^i = \sigma_1 \circ \ldots \circ \sigma_i:X_i \to X_0$ has
only simple normal crossings and $C_i$ has simple normal crossings
with $E^i,$ and \item the total transform
$(\sigma^r)^*(\mathcal{I})$ is the ideal of a simple normal
crossings divisor $E.$ If the subscheme determined by
$\mathcal{I}$ has no components of codimension one, then $E$ is a
natural combination of the irreducible components of the divisor
$E^r.$
\end{enumerate}
\end{theorem}

\begin{remark}
In order to denote the total transform
$(\sigma^r)^*(\mathcal{I}),$ other authors may use the notation
$\mathcal{I}\mathcal{O}_X.$ If $\mathcal{I}$ has components of
codimension one, we can write the total transform as a product of
two (principal) ideals: the support of the first one is the
exceptional locus, where the support of the second one is formed
by the irreducible components of the total transform that are not
contained in the exceptional locus. This second ideal is the `weak
transform' of $\mathcal{I}.$
\end{remark}

When we have a principalization $\sigma = \sigma^r,$ we can define
numerical data $(N, \nu)$ for each component of the support of
$\sigma^*(\mathcal{I})$ such that for every $b \in X$ there exist
local coordinates $(y_1, \ldots , y_n)$ which satisfy the
following conditions:
\begin{itemize}
\item if $E_1, \ldots , E_p$ are the irreducible components of the
divisor $E$ containing $b,$ we have on some neighbourhood of $b$
that $E_i$ is given by $y_i=0$ for $i=1, \ldots , p$ , \item
$\sigma^*(\mathcal{I}) \mbox{ is generated by } \varepsilon (y)
\prod_{i=1}^py_i^{N_i},$ and \item $\sigma^*(dx_1 \wedge \ldots
\wedge dx_n) = \eta(y) \prod_{i=1}^py_i^{\nu_i-1}dy_1 \wedge
\ldots \wedge dy_n,$
\end{itemize}
where $\varepsilon(y)$ and $\eta(y)$ are units in the local ring
of $X$ at $b.$ \\

We can associate a topological zeta function to an ideal
$\textbf{f}=(f_1, \ldots, f_l),$ where we suppose that $0 \in
\mbox{Supp}(\textbf{f}).$ We use the numerical data that originate
from a chosen principalization to define the local
\emph{topological zeta function}
$$Z_{top,\textbf{f}}(s) := \sum_{I\subset T} \chi(E_I^{\circ} \cap \sigma^{-1}(0))
\prod_{i\in I} \frac1{\nu_i + s N_i},$$ with $E_i(N_i, \nu_i)$ for
$i \in T$ the components of the support of the total transform of
$\textbf{f},$ and again $E_I^{\circ} = (\cap_{i \in I} E_i)
\backslash (\cup_{j\not\in I}E_j).$

When $l=1,$ Denef and Loeser showed in \cite{DL-zeta fie onafh van
res} that the expression above does not depend on the chosen
resolution by expressing it as a limit of $p$-adic Igusa zeta
functions. They introduced later in \cite{DL-motivische
zetafuncties}, still for $l=1,$ the motivic zeta function of $f,$
which is intrinsically defined. It has however a formula of the
same kind as above in terms of a resolution. Specializing this
formula to Euler characteristics yields the topological zeta
function of $f.$

One can associate more generally a motivic zeta function to an
ideal and obtain a similar formula in terms of a principalization
using the same argument as in \cite{DL-motivische zetafuncties}.
Again specializing to Euler characteristics yields the defining
expression above for the topological zeta function of an ideal.
This generalization to ideals is mentioned in
\cite[(2.4)]{veys-zuniga}.

Alternatively, one can check that this expression is independent
of the chosen principalization by verifying that it is invariant
under a blow-up with allowed center (this is straightforward) and
then applying the Weak Factorization Theorem of W\l odarczyk et
al. \cite{AKMW-factorization}. Note that in dimension 2 one does
not need the Weak Factorization Theorem since there is a minimal
principalization.

\begin{remark}
As in the case of one polynomial, there is also a \emph{global}
version of this zeta function, where we replace $E_I^{\circ} \cap
\sigma^{-1}(0)$ by $E_I^{\circ}.$ However, in this paper we will
work with the local one.
\end{remark}

Now we can ask the same questions for the topological zeta
function of an ideal in $\mathbb{C}[x,y] $ as we mentioned for the
case of one polynomial: how can we determine which candidate poles
are poles? Which rational numbers occur as poles of a zeta
function of an ideal in dimension two?

Theorem \ref{stelling polen zetafie idealen} will answer the first
question as a generalization of the result of the second author
for the topological zeta function of a curve. It turns out that
you can determine which candidate poles are poles by drawing an
intersection diagram of the $E_i$ associated to the minimal
principalization \emph{together with their numerical data.} In the
case of one polynomial, a component of the strict transform as
well as an exceptional variety that intersects at least three
times an other component, give rise to a pole. This will still
hold for the topological zeta function of an ideal in dimension
two. But this time it is not true that an exceptional variety that
intersects once or twice an other component never causes a pole.
Sometimes it will, sometimes it won't. To solve this question, we
will associate a ``generic" curve to the ideal and we will prove
that a principalization of the ideal also gives an embedded
resolution of this curve with the same numerical data. Afterwards,
we show how these numerical data tell us whether or not a
candidate pole is a pole in this case.

Further on in this paper we will answer the second question. We
will show that the possible poles of a zeta function of an ideal
in dimension two, are exactly the rational numbers in $[-1,0) \cup
\{-1-\frac1i \, | \, i \in \mathbb{Z}_{>0}\}$ (see Theorem
\ref{stelling mogelijke polen}).

In the end, we will also draw conclusions about poles of other
zeta functions of ideals in dimension two. In fact, we can say
that the same results as we prove for the topological zeta
function, are also true for the Hodge and the motivic zeta
function and for most $p$-adic Igusa zeta functions. We don't need
to prove these statements separately, but we can extract them out
of the results for the topological zeta function.

\section{Resolution of a generic curve}

Let $\textbf{f} = (f_1, \ldots , f_l)$ be an ideal in
$\mathbb{C}[x,y].$ We suppose in this section that $l>1.$ Then we
can look at the linear system $\{\lambda_1f_1 + \ldots +
\lambda_lf_l \, | \, \lambda_i \in \mathbb{C} \mbox{ for } i =1,
\ldots,l \}.$ A generic curve of $\textbf{f}$ is a general element
of this linear system. So actually, the definition of a generic
curve of an ideal is dependent on the generators we use to
represent the ideal.

\begin{lemma}
A series of blow-ups used to principalize an ideal of
$\mathbb{C}[x,y]$, also gives an embedded resolution of a generic
curve of this ideal.
\end{lemma}

\begin{remark}
This resolution will -in general- not be minimal, but we can still
use a lot of the results about the numerical data of an embedded
resolution and use them in our context.
\end{remark}

\noindent \emph{Proof.} When we start with an ideal
$\mathcal{I}=(f_1, \ldots , f_l) \subset \mathbb{C}[x,y],$ we can
first determine whether there are common components among the
$f_i$ and put them together. So we will write $$\mathcal{I} =
(h)(f_1', \ldots , f_l')$$ with $(f_1', \ldots , f_l')$ a finitely
supported ideal.

We need two chains of blow-ups to have a principalization:
\begin{description}
\item{(A)} a composition of blow-ups $\sigma : \tilde{X} \to
\mathbb{C}^2$ to transform $(f_1', \ldots , f_l')$ in a locally
principal ideal, and \item{(B)} a series of blow-ups $\tau : X \to
\tilde{X}$ to desingularize the strict transform of $h=0$ and make
it have normal crossings with all exceptional curves.
\end{description}

We will look now at the situation after the first series of
blow-ups. The ideal $\sigma^*\mathcal{I} = (f_1^*, \ldots ,
f_l^*)$ is locally principal. So in every point $b \in \tilde{X}$
we have local coordinates $(y_1, \ldots , y_n)$ and a generator
$g(y)$ such that $$f_i^*(y) = g(y) \tilde{f_i}(y)$$ for $i=1,
\ldots , l.$ Moreover, we know that there exist regular functions
$\mu_i(y)$ on $\tilde{X}$  to write that $g(y) = \sum_{i=1}^l
\mu_i(y) f_i^*(y).$ So $g(y) = g(y)\sum_{i=1}^l \mu_i(y)
\tilde{f_i}(y)$ and $1 = \sum_{i=1}^l \mu_i(y) \tilde{f_i}(y).$ We
can conclude that the $\tilde{f_i}(y)$ don't have a common zero.

We study the linear system $\{\lambda_1 \tilde{f}_1 + \ldots +
\lambda_l \tilde{f}_l =0 \, | \, \lambda_i \in \mathbb{C}\mbox{
for } i=1,\ldots , l \}.$ This is a linear system without base
points. By Bertini's theorem (see e.g.  \cite[Theorem
6.10]{Jouanolou-Bertini}) we know that a general element of the
system is non-singular and connected.

We can also restrict the linear system to an exceptional curve or
to a component of the strict transform of $h=0.$ (Note that there
are a finite number of such varieties.) On these curves, we get a
new linear system without base points. We can again use the
theorem of Bertini to say that a general element is non-singular.
In this case, this means that every intersection point of a
general element of the original linear system with a component of
the strict transform of $h=0$ or with an exceptional curve has
intersection multiplicity one.

Now we look at the following set of points: intersection points of
an exceptional curve with a component of a strict transform of
$h=0$ and singular points of the strict transform of $h=0.$ This
is a finite set. A general element of the linear system doesn't
contain any of them.

We use all this to conclude the following: if we take a generic
curve $\lambda_1 f_1' + \ldots + \lambda_l f_l'=0$ with $\lambda_1
, \ldots , \lambda_l \in \mathbb{C} $ (necessarily reduced by
Bertini's Theorem), we can suppose that the strict transform of
this curve after the first series of blow-ups (which is locally
given by $\lambda_1 \tilde{f}_1 + \ldots + \lambda_l \tilde{f}_l
=0$) is non-singular, intersects the strict transform of $h=0$ and
the exceptional curves transversely, and doesn't contain any of
the points in the mentioned set. This implies that after series
(B), the components of the strict transform of the generic curve
$\lambda_1 f_1 + \ldots + \lambda_l f_l =0$ are still non-singular
and the transform $( \sigma \circ \tau)^*(\lambda_1 f_1 + \ldots +
\lambda_l f_l)=0$ is a normal crossings divisor. So a
principalization of the ideal $(f_1 , \ldots , f_l)$ gives also an
(in general non-minimal) embedded resolution of a generic curve
$\lambda_1 f_1 + \ldots + \lambda_l f_l=0.$ \qed

\begin{remark}
This lemma is well-known. We stated and proved it in dimension
two, but one can do the same in higher dimensions. In our proof,
we made a separation in two series of blow-ups. This is not really
necessary and in higher dimensions one better avoids this.
However, we chose to make this break to get a clearer view on the
role of the common component(s) $h=0.$
\end{remark}

\begin{example} \label{vb princ-res}
We will study the ideal $(x^4y,x^7+xy^4) \subset \mathbb{C}[x,y].$
We take the generic curve $x^4y+x^7+xy^4$ of this ideal and we
perform the same blow-ups as are used to principalize the ideal.

\vspace{0.5cm} \noindent \begin{tabular}{l|c|c|}
 & principalization of  & embedded resolution of  \\
& $(x^4y,x^7+xy^4)$  & $x^4y+x^7+xy^4=0$ \\

 \hline \hline
 \underline{Chart 1} $(x,xy)$ & $x^5(y, x^2+y^4)$ & $x^5(y + x^2+y^4)$
 \\
  & $E_1 \leftrightarrow x = 0$ & $E_1 \leftrightarrow x = 0$ \\
  \hline
 \underline{Chart 2} $(xy,y)$ & $xy^5(x^3,x^6 y^2+1)$ & $xy^5(x^3+x^6 y^2+1)$
 \\
& $E \leftrightarrow x = 0$ & $E \leftrightarrow x = 0$ \\
  & $E_1 \leftrightarrow y = 0$ & $E_1 \leftrightarrow y = 0$ \\
  \hline

  \underline{Chart 1.1} $(x,xy)$ & $x^6(y,x+ x^3y^4)$ & $x^6(y+x+ x^3y^4)$
 \\
  & $E_2 \leftrightarrow x = 0$ & $E_2 \leftrightarrow x = 0$ \\
   \hline
 \underline{Chart 1.2} $(xy,y)$ & $x^5y^6(1, x^2y+y^3)$ & $x^5y^6(1+ x^2y+y^3)$
 \\
  & $E_1 \leftrightarrow x = 0$ & $E_1 \leftrightarrow x = 0$ \\
  & $E_2 \leftrightarrow y = 0$ & $E_2 \leftrightarrow y = 0$ \\
  \hline

 \underline{Chart 1.1.1} $(x,xy)$ & $x^7(y, 1+x^6y^4)$ & $x^7(y+ 1+x^6y^4)$
 \\
  & $E_3 \leftrightarrow x = 0$ & $E_3 \leftrightarrow x = 0$ \\
  \hline

   \underline{Chart 1.1.2} $(xy,y)$ & $x^6y^7(1, x+x^3y^6)$ & $x^6y^7(1+ x+x^3y^6)$
 \\
 & $E_2 \leftrightarrow x = 0$ & $E_2 \leftrightarrow x = 0$ \\
    & $E_3 \leftrightarrow y = 0$ & $E_3 \leftrightarrow y = 0$ \\
  \hline

\end{tabular}

\vspace{1cm}

\noindent We can also construct the intersection diagram of this
principalization and resolution, together with the numerical data
$(N, \nu).$

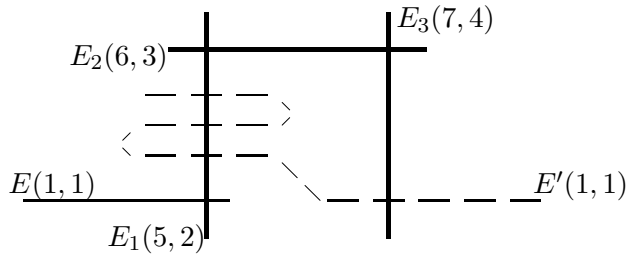
\begin{figure}[h]
\unitlength=1mm
\begin{center}
\begin{picture}(90,42)(10,0)
\linethickness{0.5mm}

\put(32,7){\line(0,1){30}}
\put(27,32){\line(1,0){34}}
\put(56,7){\line(0,1){30}}

\put(6,13){$E(1,1)$}
\put(19,6){$E_1(5,2)$}
\put(14,30){$E_2(6,3)$}
\put(57,35){$E_3(7,4)$}

\put(75,13){$E'(1,1)$}

\linethickness{0.15mm}
\put(8,12){\line(1,0){27}}

\put(24,26){\line(1,0){4}}
\put(30,26){\line(1,0){4}} \put(36,26){\line(1,0){4}}

\put(42,25){\line(1,-1){1}}
\put(42,22){\line(1,1){1}}

\put(24,22){\line(1,0){4}}
\put(30,22){\line(1,0){4}} \put(36,22){\line(1,0){4}}

\put(22,21){\line(-1,-1){1}}
\put(22,18){\line(-1,1){1}}

\put(24,18){\line(1,0){4}}
\put(30,18){\line(1,0){4}} \put(36,18){\line(1,0){4}}

\put(42,17){\line(1,-1){2}}
\put(45,14){\line(1,-1){2}}

\put(48,12){\line(1,0){4}} \put(54,12){\line(1,0){4}}
\put(60,12){\line(1,0){4}} \put(66,12){\line(1,0){4}}
\put(72,12){\line(1,0){4}}

\end{picture}
\end{center}
\end{figure}

\end{example}

\vspace{-1cm}

\noindent The curves $E$ and $E'$ are the components of the strict
transform of the generic curve. The first one is also the support
of the weak transform of the ideal, the second one does not occur
in the principalization.

\begin{remark}
\label{opm num data princ-res} In this example you can also see
that the numerical data of the principalization and those of the
resolution are the same. This is true in general. The equality of
the $\nu_i$ is obvious, the  $N_i$ are equal since for general
$\lambda_1 , \ldots , \lambda_l,$ the vanishing order of a divisor
$E$ along $\lambda_1 f_1 + \ldots + \lambda_l f_l$  is equal to
the minimum of the vanishing orders of $E$ along the $f_i.$
\end{remark}

\begin{remark} \label{opm toegelaten opblazingen}
Although the embedded resolution of the generic curve is in
general not minimal, not every blow-up is allowed in the minimal
principalization. We will only blow up with center on the
intersection of at least one exceptional curve with the support of
the weak transform of the ideal. Note that this means that
`superfluous' blowing-ups in the non-minimal embedded resolution
of our generic curve have center on the intersection of the
exceptional locus with the strict transform of the generic curve.
\end{remark}

\section{Relations between numerical data}

For the numerical data of an embedded resolution of a generic
curve of the ideal $(f_1, \ldots , f_l) \subset \mathbb{C}[x,y]$,
we know that the following relation holds: when $E(N, \nu)$ is an
exceptional curve that intersects $k$ times other components $E_i
(N_i, \nu_i)$ and $\alpha_i = \nu_i - \frac{\nu}{N}N_i$ for $i=1,
\ldots, k,$ then
$$\sum^{k}_{i=1} \alpha_i = k-2.$$
This relation between the numerical data was proved by Loeser in
\cite{Loeser num data} and generalized by the second author in
\cite{relaties num data}.

The intersection diagram with the numerical data of the
principalization is almost the same as the one that arises from
the (in general non-minimal) resolution of the generic curve $h
\cdot (\lambda_1 f_1' + \ldots + \lambda_l f_l')=0.$ Here we use
again the notation of the previous section, so we suppose that
$(f_1, \ldots , f_l)=(h)(f_1', \ldots , f_l'),$ with $(f_1',
\ldots , f_l')$ a finitely supported ideal. The only difference
between the two intersection diagrams is that the strict transform
of $\lambda_1 f_1' + \ldots + \lambda_l f_l' = 0$ disappears in
the principalization.

So we can divide the $k$ intersections of an exceptional curve of
an embedded resolution of the generic curve in two groups: there
are $n$ intersections with the strict transform of $\lambda_1 f_1'
+ \ldots + \lambda_l f_l' = 0$ and $m=k-n$ intersections that are
preserved in the intersection diagram of the principalization of
the ideal. Since we know that the first mentioned curve has
numerical data $(1,1)$, we can write -after renumbering the
intersections- that
$$\sum^{m}_{i=1} \alpha_i + n(1-\frac{\nu}{N}) = m+n-2,$$
or
\begin{equation}\label{alpha-relatie}
\sum^{m}_{i=1} \alpha_i = m-2+\frac{\nu n}{N}.
\end{equation}

\begin{proposition}
Let $E(N,\nu)$ be an exceptional curve of a principalization of
$(f_1, \ldots, f_l) \subset \mathbb{C}[x,y],$ intersecting
$E_i(N_i, \nu_i)$ for $i=1, \ldots ,m,$ and set $\alpha_i =
\nu_i-\frac{\nu}{N}N_i$ for all $i \in \{ 1, \ldots, m\}.$ Then
$-1 \leq \alpha_i < 1$ for every $i.$ Moreover, $\alpha_i = -1$
only occurs when $m=1.$
\end{proposition}

\noindent \emph{Proof.} This proposition has been proven by Loeser
in \cite[Proposition II.3.1]{Loeser num data} for the numerical
data of minimal embedded resolutions. Since we already noticed
that the numerical data of the principalization and the (possibly
non-minimal) embedded resolution of a generic curve are the same
(see Remark \ref{opm num data princ-res}), we can look at these
data as if they were coming from a resolution of the generic
curve.

We can divide the exceptional curves in two groups: the ones that
were first created are part of the minimal embedded resolution of
the generic curve. As a consequence of the mentioned theorem of
Loeser, the $\alpha_i$ that originate from these will satisfy the
condition $-1 \leq \alpha_i < 1.$ The second group of blow-ups
will have center on the intersection of one exceptional curve and
the strict transform of the generic curve. Moreover, since we
suppose that we have already an embedded resolution, we know that
the multiplicity of the generic curve in the center of the blow-up
is one.

So we only need to look at the following situation:

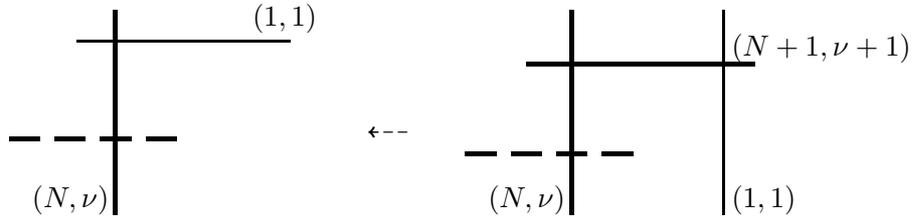
\begin{picture}(90,35)(12,3)
\unitlength=1mm \linethickness{0.5mm}

\put(13,15){\line(1,0){4}}
\put(19,15){\line(1,0){4}}
\put(25,15){\line(1,0){4}}
\put(31,15){\line(1,0){4}}
\put(27,5){\line(0,1){27}}
\linethickness{0.2mm}
\put(22,28){\line(1,0){28}}
\put(16,6){$(N,\nu)$}
\put(45,30){$(1,1)$}

\put(60,15){$\dashleftarrow$}

\linethickness{0.5mm}

\put(73,13){\line(1,0){4}}
\put(79,13){\line(1,0){4}}
\put(85,13){\line(1,0){4}}
\put(91,13){\line(1,0){4}}
\put(87,5){\line(0,1){27}}
\put(81,25){\line(1,0){30}}
\linethickness{0.2mm}
\put(107,5){\line(0,1){27}}

\put(76,6){$(N,\nu)$}
\put(108,26){$(N+1, \nu +1)$}
\put(108,6){$(1, 1)$}

\end{picture}

\noindent We can suppose that $-1 \leq 1-\frac{\nu}{N} < 1$ (or
that $ 0 < \frac{\nu }N \leq 2$) and we only need to show that

\begin{enumerate}[(i)]
\item $-1 \leq \nu +1 - \frac{\nu}N (N+1) < 1, $ \item $-1 \leq
\nu - \frac{\nu+1}{N+1} N < 1$ \hspace{0.3cm} and \item $-1 \leq 1
- \frac{\nu+1}{N+1} < 1.$
\end{enumerate}
This is straightforward. \qed

\begin{corollary} \label{gevolg voor geordende boom}
Let $E(N,\nu)$ be an exceptional curve of a principalization
$\sigma$ of an ideal $\mathcal{I} \subset \mathbb{C}[x,y].$
Suppose that $E$ intersects the other components $E_i(N_i, \nu_i)$
for $i=1, \ldots, m$ of the total transform $\sigma^*\mathcal{I}.$
Let $\alpha_i = \nu_i-\frac{\nu}{N}N_i$ for all $i \in \{ 1,
\ldots, m\}.$ Then we have the following statements.

\begin{enumerate} \item At most one $E_i, 1\leq i \leq m,$ occurs such that $\alpha_i
< 0.$ \item If $m \geq 3,$ then there is at most one $i$ such that
$\alpha_i \leq 0.$ \item If $m=2,$ we see that $\frac{\nu_1}{N_1}
< \frac{\nu}{N} \Rightarrow \frac{\nu}{N} < \frac{\nu_2}{N_2}.$
\end{enumerate}
\end{corollary}
\noindent This is a direct consequence of the previous proposition
and equation (\ref{alpha-relatie}). In \cite{veys polen top
zetafunctie} there are almost the same results for the numerical
data of an embedded resolution of a curve. However, the analogue
of the third statement in that context is an equivalence instead
of an implication. Roughly said, this is due to the presence of
the positive term $\frac{\nu n}{N}$ in our equation
(\ref{alpha-relatie}).

The mentioned corollary in \cite{veys polen top zetafunctie} is
used there to determine the `ordered tree'-structure of the
resolution graph. The same can be done in our case. We can draw a
dual principalization graph by associating a vertex to every
exceptional curve and every (analytically irreducible) component
of the support of the weak transform. For each intersection we
have an edge, connecting the corresponding vertices.

By using Corollary \ref{gevolg voor geordende boom}, it is not so
difficult to derive the next proposition. For example, this can be
done as in \cite[Theorem 3.3]{veys polen top zetafunctie}.
\begin{proposition}
\label{prop ordered tree} The part of the dual principalization
graph where $\frac{\nu}N$ is minimal, is connected. Moreover, when
we follow a path that moves away from this minimal part, the ratio
$\frac{\nu}N$ will strictly increase.
\end{proposition}

\section{Poles of a zeta function of an ideal}
In this section we always consider ideals $\mathcal{I} \subset
\mathbb{C}[x,y]$ with $0 \in \mbox{Supp}(\mathcal{I}).$ Since we
study the local topological zeta function associated to
$\mathcal{I},$ we need in fact only a principalization of
$\mathcal{I}$ in the neighbourhood of $0.$

We know that the only possible poles of the topological zeta
functions are rational numbers $-\frac{\nu}N$ with $(N,\nu)$
numerical data of components of the minimal principalization. We
can see that the largest candidate pole plays a special role. The
following arguments show that it is always a pole. If there are
different components with this maximal ratio $-\frac{\nu}N,$ these
components need to intersect and we find a pole of order two.
Moreover, this is the only value where a pole of order two is
possible. This is a consequence of the `ordered tree'-structure of
the graph (see Proposition \ref{prop ordered tree}). When there is
only one component $E(N, \nu)$ with this minimal ratio, we have a
candidate pole of order one. Its residue is then given by
$$R = \frac1N\left(2-m+\sum_{i=1}^m\frac1{\alpha_i}\right),$$
where we suppose that $E$ intersects $m$ times other components
$E_i(N_i,\nu_i)$ $(i=1,\ldots,m)$ of the principalization, and
$\alpha_i = \nu_i - \frac{\nu}N N_i.$ When $\frac{\nu}N$ is
minimal, then $0<\alpha_i <1$ for every $i,$ so $R>0$ and $-\frac{\nu}N$ is a pole.\\
\\
Not every other candidate pole gives rise to a pole. For the
topological zeta function associated to a curve in $\mathbb{C}^2,$
the second author proved the following theorem in \cite{veys polen
top zetafunctie}.
\begin{theorem} \label{stelling Veys 3 keer snijden} Let $f \in
\mathbb{C}[x,y]$ be a non-constant polynomial satisfying $f(0)=0,$
and let $h:X \to \mathbb{C}^2$ be the minimal embedded resolution
of $f^{-1}\{0\}$ in a neighbourhood of 0. Let $E_i(N_i, \nu_i)$ be
the irreducible components of $h^{-1}(f^{-1}\{0\})$ with their
associated numerical data. We have that $s_0$ is a pole of
$Z_{top,f}(s)$ if and only if $s_0 = -\frac{\nu_i}{N_i}$ for some
exceptional curve $E_i$ intersecting at least three times other
components or $s_0 = -\frac1{N_i}$ for some irreducible component
$E_i$ of the strict transform of $f=0.$
\end{theorem}
This gives a criterion to filter the poles out of the series of
candidate poles. The next theorem will do the same for the
topological zeta function associated to an ideal in
$\mathbb{C}[x,y].$ With this theorem we can easily determine the
poles of the zeta function when we have the principalization of
the ideal.

\begin{theorem}
\label{stelling polen zetafie idealen} Let $\mathcal{I} \subset
\mathbb{C}[x,y]$ be an ideal satisfying $0 \in
\mbox{Supp}\,(\mathcal{I})$ and $\sigma: X \to \mathbb{C}^2$ the
minimal principalization of $\mathcal{I}$ in a neighbourhood of 0.
Let $E_\bullet(N_\bullet,\nu_\bullet)$ be the components of the
support of the total transform $\sigma^*\mathcal{I}$ with their
associated numerical data.

The rational number $s_0$ is a pole of the local topological zeta
function of $\mathcal{I}$ if and only if one of the following
conditions is satisfied:
\begin{enumerate}
\item $s_0 = -\frac1{N}$ for a component $E(N, \nu)$ of the
support of the weak transform of $\mathcal{I}$; \item $s_0 =
-\frac{\nu}{N}$ for $E(N, \nu)$ an exceptional curve that
intersects no other component; \item $s_0 = -\frac{\nu}{N}$ for
$E(N, \nu)$ an exceptional curve that intersects once another
component $E_i(N_i, \nu_i)$ with $\nu_i - \frac{\nu}{N} N_i \neq
-1;$ \item $s_0 = -\frac{\nu}{N}$ for $E(N, \nu)$ an exceptional
curve that intersects two times other components $E_i(N_i, \nu_i)$
and $E_j(N_j, \nu_j)$ with $(\nu_i - \frac{\nu}{N} N_i) + (\nu_j -
\frac{\nu}{N} N_j)  \neq 0;$ \item $s_0 = -\frac{\nu}{N}$ for
$E(N, \nu)$ an exceptional curve that intersects at least three
times other components.
\end{enumerate}
\end{theorem}

\begin{remark}
In the proof we will work with the following notation. If $E(N,
\nu)$ is a curve in the support of the total transform of
$\mathcal{I}$ that intersects once another curve $E_i(N_i,
\nu_i),$ we write $\alpha = \nu_i - \frac{\nu}{N}N_i.$ If $E(N,
\nu)$ intersects the curves $E_{i_1}(N_{i_1}, \nu_{i_1}),
E_{i_2}(N_{i_2}, \nu_{i_2}),$ $\ldots ,$ $E_{i_m}(N_{i_m},
\nu_{i_m}),$ we write $\alpha_j = \nu_{i_j} -
\frac{\nu}{N}N_{i_j}.$
\end{remark}

\noindent \emph{Proof.} We have already said that the only
possible pole of order two is the largest candidate pole. We can
see that if $s_0$ is maximal, at least one of the five conditions
is satisfied.

Now we will calculate the contribution to the residue of $s_0$ as
a pole of order one in the various cases. (We can suppose that all
$\alpha_i \neq 0.$) We will see that the five situations of the
theorem are the only ones where that contribution is non-zero.
Moreover we will show that this contribution is negative, unless
$s_0$ is the largest candidate pole. Notice that this last
condition corresponds with ``every $\alpha_i > 0".$

Suppose that $s_0 = -\frac1{N}$ for a component $E(N, \nu)$ of the
support of the weak transform. Such a component only intersects
one exceptional curve $E_i(N_i, \nu_i)$ and we see that the
contribution to the residue of a pole of order one is $R=
\frac{1}{N\alpha}.$ This is positive if $s_0$ is the largest pole
and negative otherwise.

If we have an exceptional curve $E(N, \nu)$ that doesn't intersect
any other component, we know that this is the only curve in the
principalization. So the topological zeta function is given by
$\frac{2}{\nu + s N}$ and the value $s_0 = -\frac{\nu}{N}$ is a
pole.

Let $s_0 = -\frac{\nu}{N}$ for $E(N, \nu)$ an exceptional curve
that intersects once another component. The contribution to the
residue for a pole of order one, is $R=
\frac1N(2-1+\frac1{\alpha}).$ So we see immediately that $$R=0
\Leftrightarrow \alpha = -1.$$ The case $\alpha = 0$ is excluded,
because we only look at candidate poles of order one. If $\alpha >
0,$ then $R>0.$ If $\alpha < 0$ and $\alpha \neq -1,$ we can use
$\alpha > -1$ to conclude that $R<0.$

Now suppose that $s_0 = -\frac{\nu}{N}$ for $E(N, \nu)$ an
exceptional curve that intersects two times other components. In
this case $R = \frac1N(2-2+\frac1{\alpha_1} + \frac1{\alpha_2}
)=\frac1N \frac{\alpha_1+\alpha_2}{\alpha_1\alpha_2}.$ We conclude
that $$R = 0 \Leftrightarrow \alpha_1+\alpha_2=0.$$ If
$\alpha_1+\alpha_2\neq 0,$ we can use equation
(\ref{alpha-relatie}) to know that $\alpha_1+\alpha_2 > 0.$ So we
are only interested in the sign of $\alpha_1\alpha_2$ to know the
sign of $R.$ We know that $\alpha_1$ and $\alpha_2$ can't be both
negative. If they are both positive, then $s_0$ is the largest
pole and $R>0.$ In the other case we have $R<0.$

The next case is where $s_0 = -\frac{\nu}{N}$ for $E(N, \nu)$ an
exceptional curve that intersects at least three times an other
component. Here, the contribution to the residue $R =
\frac1N(2-m+\sum_{i=1}^m\frac1{\alpha_i})$ is always non-zero. If
every $\alpha_i >0,$ we can easily conclude that $R>0.$ When there
is a $\alpha_i < 0,$ we can use the results for the resolution of
a curve that are written in \cite[Proposition 2.8]{veys polen top
zetafunctie} to see that $R' := \frac1N
\left(2-(m+n)+\sum_{i=1}^m\frac1{\alpha_i} +
\frac{n}{1-\frac{\nu}{N}} \right)<0.$ (Here we use the notation of
the previous section.) Because we know that there can exist at
most one negative $\alpha,$ we can deduce that $0 < 1 -
\frac{\nu}{N}< 1$ and $$R = R' +\frac1N \left( n-
\frac{n}{1-\frac{\nu}{N}}\right) < 0.$$

From these calculations we can conclude that all contributions to
the residue are negative if $s_0$ is not maximal. So if one of
these contributions is non-zero, the total residue is non-zero and
$s_0$ is a pole of order one. \qed

\begin{remark}
When we work with one element $f \in \mathbb{C}[x,y]$ instead of
an ideal, only the first and the last case can occur. Moreover, a
principalization of the ideal $(f)$ is the same as an embedded
resolution of the curve given by $f=0.$ So this theorem is a
generalization of Theorem \ref{stelling Veys 3 keer snijden}.
\end{remark}

\begin{example}
We will continue Example \ref{vb princ-res} and calculate
explicitly the topological zeta function of the ideal $\textbf{f}
= (x^4y,x^7+xy^4) \subset \mathbb{C}[x,y].$ With the calculations
done in the previous example, we can see that
$$Z_{top, \textbf{f}}(s) = \frac1{4+7s} + \frac1{(3+6s)(4+7s)} +
\frac1{(2+5s)(3+6s)} + \frac1{(1+s)(2+5s)}.$$  With a little
calculation, we can simplify this expression to
$$Z_{top, \textbf{f}}(s) = \frac{5s^2+16s+8}{(4+7s)(2+5s)(1+s)},$$
which implies that the poles of this function are $-4/7, -2/5$ and
$-1.$

We can obtain the same result by using Theorem \ref{stelling polen
zetafie idealen} in the following way:
\begin{itemize}
\item $E(1,1)$ is a component of the support of the weak
transform, so $-1$ is a pole; \item $E_1(5,2)$ intersects twice
other components, with $(1-\frac25 ) + (3-\frac25 6) \neq 0,$ so
$-2/5$ is a pole; \item $E_2(6,3)$ also has two intersections with
other components, this time with $(4-\frac367)+(2-\frac365) = 0,$
hence the candidate pole $-1/2$ is no pole; \item $E_3(7,4)$
intersects one other component with $3-\frac476 \neq -1,$ so this
gives the last pole $-4/7.$
\end{itemize}
\end{example}

\section{Determination of all possible poles}
In this section, we will determine which numbers can occur as a
pole of a topological zeta function associated to an ideal in
dimension 2. For the topological zeta function of a curve, this
question has been answered in \cite{SV smallest poles} and
\cite{artikel Ann & co}.

In the first article, Segers and the second author proved that the
poles smaller than $-\frac12$ are given by $\{-\frac12 - \frac1i
\, | \, i \in \mathbb{Z}_{>1} \}.$ In the second article,
Lemahieu, Segers and the second author showed that every rational
number in the interval $[-\frac12,0)$ is a pole of a zeta function
of a curve. This determines all possible poles.

We will prove an analogue of these results for the topological
zeta function of an ideal in dimension 2.
\\

For an exceptional variety $E(N, \nu)$ of the minimal embedded
resolution of a curve, one can show that $\nu \leq N.$
Analogously, we prove the following proposition.

\begin{proposition}
Let $E(N, \nu)$ be an exceptional curve of the minimal
principalization of an ideal in $\mathbb{C}[x,y].$ Then the
numerical data satisfy
$$\nu \leq N + 1.$$
\end{proposition}

\noindent \emph{Proof.} We will prove this proposition by
induction. If $E_1(N_1, \nu_1)$ is the first created exceptional
curve, then $\nu_1=2$ and $N_1 \geq 1,$ so the statement is
proven. Now we will suppose that the inequality is satisfied for
all already created exceptional curves and we will prove it for
the next one.
\begin{itemize}
\item First, suppose that the center of the blow-up is contained
in two exceptional curves $E_{i_1}(N_{i_1}, \nu_{i_1})$ and
$E_{i_2}(N_{i_2}, \nu_{i_2}).$ Then we see that $\nu_i = \nu_{i_1}
+ \nu_{i_2}$ and $N_i = N_{i_1} + N_{i_2} + $ (minimal
multiplicity of the generators of the ideal in the center). If we
use the induction hypothesis, we see that $\nu_i \leq N_{i_1} +
N_{i_2} + 2,$ but we also know that $N_i \geq N_{i_1} + N_{i_2} +
1,$ so $\nu_i \leq N_i + 1.$ \item When only $E_{i_1}$ exists,
then $\nu_i = \nu_{i_1} + 1 \leq N_{i_1} + 2$ and $N_i = N_{i_1} +
$ (minimal multiplicity of the generators of the ideal in the
center) $ \geq N_{i_1} + 1,$ so $\nu_i \leq N_i + 1.$ \qed
\end{itemize}

\begin{remark} \label{opm beschrijving mogelijke polen }
This implies that all candidate poles of the topological zeta
function associated to an ideal are rational elements of $[-1,0)
\cup \{-1-\frac1i \, | \, i \in \mathbb{Z}_{>0}\}.$ The next
proposition will show that every rational number in this range
really occurs as a pole of a certain topological zeta function.
Hence we have a complete description of the possible poles.
\end{remark}

\begin{theorem} \label{stelling mogelijke polen}
The set  of rational numbers $s_0$ for which there exists an ideal
$\mathcal{I}\subset \C[x,y]$ such that $Z_{top, \mathcal{I}}\,(s)$
has a pole in $s_0,$ is given by $\mathbb{Q} \, \cap \, ([-1,0)
\cup \{-1-\frac1i \, | \, i \in \mathbb{Z}_{>0}\}).$
\end{theorem}

\noindent \emph{Proof.} Choose $a,b \in \mathbb{Z}_{\geq 0}$ with
$a>b.$ Look at the ideal $$(x^by\, , \, x^a+y^{b+1}) \subset
\mathbb{C}[x,y].$$ After principalization, we find the following
numerical data and intersection diagram:
$$E_1(b+1,2) \ , \ E_2(b+2,3) \ , \, \ldots \, , \
E_{a-b-1}(a-1,a-b)\ , \ E_{a-b}(a, a-b+1).$$

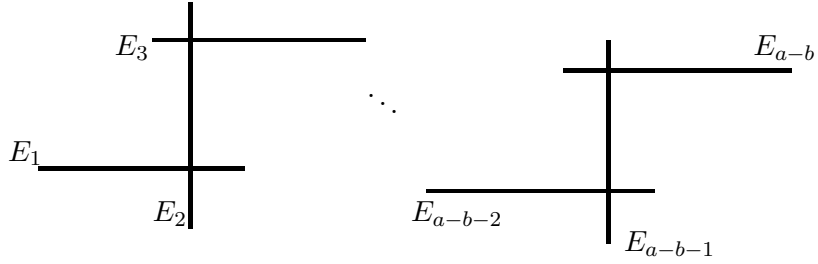
\vspace{0.5cm}

\begin{picture}(90,35)(0,0)
\unitlength=1mm
\linethickness{0.5mm}

\put(12,15){\line(1,0){27}}
\put(32,7){\line(0,1){30}}
\put(27,32){\line(1,0){28}}
\put(8,16){$E_1$} \put(27,8){$E_2$} \put(22,30){$E_3$}

\put(55,22){$\ddots$}

\put(63,12){\line(1,0){30}}
\put(87,5){\line(0,1){27}}
\put(81,28){\line(1,0){30}}

\put(61,8){$E_{a-b-2}$}
\put(89,4){$E_{a-b-1}$}
\put(106,30){$E_{a-b}$}

\end{picture}

The last exceptional variety only intersects  $E_{a-b-1}$ with
$\alpha \neq -1,$ so this causes a pole in $-\frac{a-b+1}{a}.$
Easy calculations show that this implies that every element of $\,
\mathbb{Q} \, \cap \, ([-1,0) \cup \{-1-\frac1i \, | \, i \in
\mathbb{Z}_{>0}\})$ occurs as a pole of the topological zeta
function of an ideal in $\mathbb{C}[x,y].$ \qed

\section{Other zeta functions}

There are finer variants of the topological zeta function of an
ideal. For instance, using the notation of the introduction, there
is the (local) \emph{Hodge zeta function}
$$Z_{Hod, \textbf{f}}(s)=\sum_{I \subset T} H(E_I^\circ \cap
\sigma^{-1}(0); u,v) \prod_{i \in I}
\frac{uv-1}{(uv)^{\nu_i+sN_i}-1} \in \mathbb{Q}(u,v)((uv)^{-s})$$
for the ideal $\textbf{f},$ where $H(\, \cdot \, ;u,v) \in
\mathbb{Z}[u,v]$ denotes the Hodge polynomial. Even finer is the
(local) motivic zeta function of $\textbf{f},$ which was already
mentioned in the introduction. Its formula involves classes in the
Grothendieck ring of algebraic varieties, instead of Euler
characteristics or Hodge polynomials. We refer to e.g.
\cite{DL-motivische zetafuncties}, \cite{Rodrigues} or \cite{V arc
spaces etc} for these zeta functions and their global versions
associated to one polynomial and to \cite{veys-zuniga} for ideals.

The results in this paper on poles of the topological zeta
function of an ideal in dimension 2, i.e. Theorems \ref{stelling
polen zetafie idealen} and \ref{stelling mogelijke polen}, are
also valid for the Hodge and the motivic zeta function. We chose
not to give the details here about these zeta functions, since for
results of the kind we proved, the version for the topological
zeta function is the strongest, and implies the same results for
the finer zeta functions.

The point is that the motivic zeta function specializes to the
Hodge zeta function, which in turn specializes to the topological
zeta function. (Note for instance that $H(\, \cdot \, ;1,1) =
\chi(\cdot).$) In particular, a pole of the topological zeta
function will induce a pole of the other two. (The converse is not
clear.) We refer to \cite{Rodrigues} and \cite{Rodrigues-Veys} for
the precise description of the notion of a pole for the Hodge and
the motivic zeta function. Here we should note that the analogue
of Theorem \ref{stelling polen zetafie idealen} for the finer zeta
functions also requires the verification of the following in the
context of for instance Hodge polynomials. Exceptional curves
intersecting once or twice other components such that $\alpha=-1$
or $\alpha_1+\alpha_2 = 0,$ respectively, should not contribute to
the residue of the induced candidate pole. Now this is as
straightforward as with Euler characteristics (and well known). \\

Theorems \ref{stelling polen zetafie idealen} and \ref{stelling
mogelijke polen} are also valid for (most) $p$-adic Igusa zeta
functions. We briefly introduce the necessary notation to
introduce these zeta functions and to state the precise result.

Let $K$ be a finite extension of the $p$-adic numbers with
valuation ring $R,$ maximal ideal $P,$ and residue field $\bar{K}
= R/P (\cong \mathbb{F}_q).$ Denote for $z \in K$ by $|z|$ its
standard absolute value, and put $\|z\| := \max_{1\leq i \leq l}
|z_i|$ for $z=(z_1, \ldots , z_l) \in K^l.$ Let $f_1, \ldots ,
f_l$ be polynomials in $K[x_1, \ldots , x_n].$ The (local)
\emph{$p$-adic Igusa zeta function} associated to the mapping
$\textbf{f} = (f_1 , \ldots , f_l): K^n \to K^l$ is
$$Z_{K, \textbf{f}}(s) := \int_{P^n} \| \textbf{f}(x) \|^s \, |dx|$$
for $s \in \mathbb{C}$ with $\Re(s) > 0,$ where $|dx|$ is the
usual Haar measure on $K^n.$ A global version consists in
replacing $P^n$ by $R^n.$ This function is analytic in $s$ and
admits a meromorphic continuation to $\mathbb{C}$ as a rational
function of $q^{-s}.$ This was first proved by Igusa for $l=1$
(see \cite{Igusa}). For arbitrary $l$ there are different proofs
in \cite{Meuser}, \cite{Denef - rationaliteit Poincare} and
\cite{veys-zuniga}.

Considering polynomials $f_1, \ldots , f_l$ over a number field
$F,$ one can study $Z_{K,F}(s)$ for all (non-archimedean)
completions $K$ of $F.$ For all but finitely many completions $K$
there is a concrete formula for $Z_{K,F}(s)$ in terms of a
principalization of the ideal $(f_1, \ldots , f_l),$ similar to
the formulas for the other zeta functions in this paper. This was
proved for $l=1$ by Denef in \cite{Denef - igusa zetafunctie}, and
can be generalized to arbitrary $l,$ see
\cite[(2.3)]{veys-zuniga}. (In fact the motivic zeta function
specializes to `almost all' $p$-adic zeta functions, see
\cite[(2.4)]{DL-motivische zetafuncties}.)

When the number field $F$ is large enough, then for all but
finitely many completions $K$ of $F$ we have that the analogues of
Theorems \ref{stelling polen zetafie idealen} and \ref{stelling
mogelijke polen} are valid for $Z_{K,F}(s),$ replacing `pole' by
`real pole'. One can derive this from the results for the
topological zeta function, or by completely analogous proofs.
(Previous such results for $l=1$ are in \cite{V polen igusa
krommen} and \cite{Segers - kleinste polen}.)

\footnotesize{
 
}

\end{document}